\documentclass[psamsfonts,english]{amsart}

    \usepackage[T1]{fontenc}
    \usepackage{babel}
    \usepackage{natbib}
    \usepackage{amsmath}
    \usepackage{amssymb}
    \usepackage{amsthm}
    \usepackage[mathcal]{eucal}
    \usepackage{url}

    \newtheorem{teo}{Theorem}
    \newtheorem{lem}[teo]{Lemma}
    \theoremstyle{definition}
    \newtheorem{defn}[teo]{Definition}
    \theoremstyle{remark}
    \newtheorem{rem}[teo]{Remark}

    \newcommand{\FF}{\mathbb{F}}
    \newcommand{\QQ}{\mathbb{Q}}
    \newcommand{\ZZ}{\mathbb{Z}}
    \newcommand{\RR}{\mathbb{R}}
    \newcommand{\PP}{\mathbb{P}}
    \newcommand{\jac}[1]{\mathcal{J}_{#1}}
    \newcommand{\heltal}[1]{\mathfrak{O}_{#1}}

    \DeclareMathOperator{\End}{End}
    \DeclareMathOperator{\Div}{Div}
    \DeclareMathOperator{\gal}{Gal}

\begin{document}

\title[$p$-torsion of Genus Two Curves Over $\FF_p$] 
{$p$-torsion of Genus Two Curves Over Prime Fields of Characteristic $p$}

\author[C.R. Ravnshøj]{Christian Robenhagen Ravnshøj}

\address{Department of Mathematical Sciences \\
University of Aarhus \\
Ny Munkegade \\
Building 1530 \\
DK-8000 Aarhus C}

\email{cr@imf.au.dk}

\thanks{Research supported in part by a PhD grant from CRYPTOMAThIC}

\keywords{Jacobians, hyperelliptic curves, complex multiplication, cryptography}

\subjclass[2000]{Primary 14H40; Secondary 11G15, 14Q05, 94A60}


\begin{abstract}
Consider the Jacobian of a hyperelliptic genus two curve defined over a prime field of characteristic $p$ and with
complex multiplication. In this paper we show that the $p$-Sylow subgroup of the Jacobian is either trivial or of
order~$p$.
\end{abstract}

\maketitle

\section{Introduction}

In elliptic curve cryptography it is essential to know the number of points on the curve. Cryptographically we are
interested in elliptic curves with large cyclic subgroups. Such elliptic curves can be constructed. The construction is
based on the theory of complex multiplication, studied in detail by \cite{atkin-morain}. It is referred to as the
\emph{CM method}.

\cite{koblitz89} suggested the use of hyperelliptic curves to provide larger group orders. Therefore constructions of
hyperelliptic curves are interesting. The CM method for elliptic curves has been generalized to hyperelliptic curves of
genus~two by \cite{spallek}, and efficient algorithms have been proposed by \cite{weng03} and \cite{gaudry}.

Both algorithms take as input a primitive, quartic CM field $K$ (see section~\ref{sec:CMfields} for the definition of a
CM~field), and give as output a hyperelliptic genus~two curve $C$ defined over a prime field $\FF_p$. A prime number
$p$ is chosen such that $p=x\overline x$ for a number $x\in\heltal{K}$, where $\heltal{K}$ is the ring of integers of
$K$. We have $K=\QQ(\eta)$ and $K\cap\RR=\QQ(\sqrt{D})$, where $\eta=i\sqrt{a+b\xi}$ and
    $$
    \xi=\begin{cases}
    \frac{1+\sqrt{D}}{2}, & \textrm{if $D\equiv 1\mod{4}$,} \\
    \sqrt{D}, & \textrm{if $D\equiv 2,3\mod{4}$}.
    \end{cases}
    $$
In this paper, the following theorem is established.

\begin{teo}\label{teo:maal}
Let $C$ be a hyperelliptic curve of genus two defined over a prime field~$\FF_p$. Assume that
$\End(C)\simeq\heltal{K}$, where $K$ is a primitive, quartic CM field as defined in
definition~\ref{def:CMfieldPrimitive}, and that the $p$-power Frobenius under this isomorphism is given by a number
in $\heltal{K_0}+\eta\heltal{K_0}$, where $\eta$ is given as above. Then the $p$-Sylow subgroup of
$\jac{C}(\FF_p)$ is either trivial or of order~$p$.
\end{teo}

\section{Hyperelliptic curves}

A hyperelliptic curve is a smooth, projective curve $C\subseteq\PP^n$ of genus at least two with a separable, degree
two morphism $\phi:C\to\PP^1$. Let $C$ be a hyperelliptic curve of genus two defined over a prime field $\FF_p$ of
characteristic~$p>2$. By the Riemann-Roch theorem there exists an embedding $\psi:C\to\PP^2$, mapping $C$ to a curve
given by an equation of the form
    $$y^2=f(x),$$
where $f\in\FF_p[x]$ is of degree six and have no multiple roots \cite[see][chapter~1]{cassels}.

The set of principal divisors $\mathcal{P}(C)$ on $C$ constitutes a subgroup of the degree 0 divisors $\Div_0(C)$. The
Jacobian $\jac{C}$ of $C$ is defined as the quotient
    $$\jac{C}=\Div_0(C)/\mathcal{P}(C).$$
Since $C$ is defined over $\FF_p$, the mapping $(x,y)\mapsto (x^p,y^p)$ is a morphism on~$C$. This morphism induces the
$p$-power Frobenius endo\-morphism $\varphi$ on the Jacobian $\jac{C}$. The characteristic polynomial $P(X)$ of
$\varphi$ is of degree four \cite[Theorem~2, p.~140]{tate}, and by the definition of $P(X)$
\cite[see][pp.~109--110]{lang59},
    $$|\jac{C}(\FF_p)|=P(1),$$
i.e. the number of $\FF_p$-rational points on the Jacobian is determined by $P(X)$.

\section{CM fields}\label{sec:CMfields}

An elliptic curve $E$ with $\ZZ\neq\End(E)$ is said to have \emph{complex multiplication}. Let $K$ be an ima\-ginary,
quadratic number field with ring of integers $\heltal{K}$. $K$ is a \emph{CM field}, and if
\mbox{$\End(E)\simeq\heltal{K}$}, then $E$ is said to have \emph{CM by $\heltal{K}$}. More generally a CM field is
defined as follows.

\begin{defn}[CM field]
A number field $K$ is a CM field, if $K$ is a totally imaginary, quadratic extension of a totally real number field
$K_0$.
\end{defn}

In this paper only CM fields of degree $[K:\QQ]=4$ are considered. Such a field is called a \emph{quartic} CM field.

\begin{rem}\label{rem:quarticCM}
Consider a quartic CM field $K$. Let $K_0=K\cap\RR$ be the real subfield of $K$. Then $K_0$ is a real, quadratic number
field, $K_0=\QQ(\sqrt{D})$. By a basic result on quadratic number fields, the ring of integers of $K_0$ is given by
$\heltal{K_0}=\ZZ+\xi\ZZ$, where
    $$
    \xi=\begin{cases}
    \frac{1+\sqrt{D}}{2}, & \textrm{if $D\equiv 1\mod{4}$,} \\
    \sqrt{D}, & \textrm{if $D\equiv 2,3\mod{4}$}.
    \end{cases}
    $$
Since $K$ is a totally imaginary, quadratic extension of $K_0$, a number $\eta\in K$ exists, such that $K=K_0(\eta)$,
$\eta^2\in K_0$. The number $\eta$ is totally imaginary, and we may assume that $\eta=i\eta_0$, $\eta_0\in\RR$.
Furthermore we may assume that $\eta^2\in\heltal{K_0}$; so $\eta=i\sqrt{a+b\xi}$, where $a,b\in\ZZ$.
\end{rem}

Let $C$ be a hyperelliptic curve of genus two. Then $C$ is said to have CM by~$\heltal{K}$, if
$\End(C)\simeq\heltal{K}$. The structure of $K$ determines whether $C$ is irreducible. More precisely, the following
theorem holds.

\begin{teo}\label{teo:reducibel}
Let $C$ be a hyperelliptic curve of genus two with $\End(C)\simeq\heltal{K}$, where $K$ is a quartic CM field. Then $C$
is reducible if, and only if, $K/\QQ$ is Galois with Galois group $\gal(K/\QQ)\simeq\ZZ/2\ZZ\times\ZZ/2\ZZ$.
\end{teo}

\begin{proof}
\cite[Proposition~26, p.~61]{shi}.
\end{proof}

{\samepage Theorem~\ref{teo:reducibel} motivates the following definition.

\begin{defn}[Primitive, quartic CM field]\label{def:CMfieldPrimitive}
A quartic CM field $K$ is called primitive if either $K/\QQ$ is not Galois, or $K/\QQ$ is Galois with cyclic Galois
group.
\end{defn}
}

The CM method for constructing curves of genus~two with prescribed endomorphism ring is described in detail by
\cite{weng03} and \cite{gaudry}. In short, the CM method is based on the construction of the class polynomials of a
primitive, quartic CM field $K$ with real subfield $K_0$ of class number $h(K_0)=1$. The prime number $p$ has to be
chosen such that $p=x\overline x$ for a number $x\in\heltal{K}$. By \cite{weng03} we may assume that
$x\in\heltal{K_0}+\eta\heltal{K_0}$.

\section{The $p$-Sylow subgroup of $\jac{C}(\FF_p)$}\label{sec:properties}

Let $K$ be a primitive, quartic CM field with real subfield $K_0=\QQ(\sqrt{D})$ of class number $h(K_0)=1$. Cf.
Remark~\ref{rem:quarticCM} we may write $K=\QQ(\eta)$, where $\eta=i\sqrt{a+b\xi}$ and
    $$\xi=\begin{cases}
    \frac{1+\sqrt{D}}{2}, & \textrm{if $D\equiv 1\mod{4}$,} \\
    \sqrt{D}, & \textrm{if $D\equiv 2,3\mod{4}$}.
    \end{cases}
    $$
Let $p$ be a prime number such that $p=x\overline x$ for a number $x\in\heltal{K_0}+\eta\heltal{K_0}$. Let $C$ be a
hyperelliptic curve of genus two defined over $\FF_p$ with $\End(C)\simeq\heltal{K}$.  Assume that the $p$-power
Frobenius under this isomorphism is given by the number
    \begin{equation}\label{eq:omega}
    \omega=c_1+c_2\xi+(c_3+c_4\xi)\eta,\quad c_i\in\ZZ.
    \end{equation}
Since the $p$-power Frobenius is of degree $p$, we know that $\omega\overline\omega=p$.

\begin{rem}\label{rem:deg4}
If $c_2=0$ in~\eqref{eq:omega}, then $\gal(K/\QQ)\simeq\ZZ/2\ZZ\times\ZZ/2\ZZ$, and $K$ is not primitive. So $c_2\neq 0$.
\end{rem}

The characteristic polynomial $P(X)$ of the Frobenius is given by
    $$
    P(X) = \prod_{i=1}^4(X-\omega_i),
    $$
where $\omega_i$ are the conjugates of $\omega$. Since the conjugates of $\omega$ are given by $\omega_1=\omega$,
$\omega_2=\overline{\omega}_1$, $\omega_3$ and $\omega_4=\overline{\omega}_3$, where
$\omega_3=c_1+c_2\xi'+(c_3+c_4\xi')\eta'$, $\eta'=i\sqrt{a+b\xi'}$ and
    $$
    \xi' =
     \begin{cases}
      -\sqrt{D}, & \textrm{ if } D\equiv 2,3 \mod{4} \\
      \frac{1-\sqrt{D}}{2}, & \textrm{ if } D\equiv 2,3 \mod{4}
     \end{cases}
    $$
it follows that
    \begin{align*}
     P(X) &= X^4-4c_1X^3+(2p+4(c_1^2-c_2^2D))X^2-4c_1pX+p^2,
     \intertext{if $D\equiv 2,3\mod{4}$, and}
     P(X) &= X^4-2cX^3+(2p+c^2-c_2^2D)X^2-2cpX+p^2,
    \end{align*}
if $D\equiv 1\mod{4}$. Here, $c=2c_1+c_2$. We notice that $4\mid P(1)=|\jac{C}(\FF_p)|$. This observation leads
to the following lemma.

\begin{lem}\label{lem:1}
Let $C$ be a hyperelliptic curve of genus two defined over a prime field~$\FF_p$ of characteristic $p>5$. Assume
that $\End(C)\simeq\heltal{K}$ and that the $p$-power Frobenius under this isomorphism is given by a number in
$\heltal{K_0}+\eta\heltal{K_0}$, where $\eta$ is given as in remark~\ref{rem:quarticCM}. Then the $p$-Sylow
subgroup of $\jac{C}(\FF_p)$ is either trivial or of order~$p$.
\end{lem}

\begin{proof}
Assume $p^2\mid N=|\jac{C}(\FF_p)|$. Since $|\omega_i|=\sqrt{p}$, we know that
    $$N=P(1)=\prod_{i=1}^4(1-\omega_i)\leq (1+\sqrt{p})^4=p^2+4p\sqrt{p}+6p+4\sqrt{p}+1.$$
Hence, $\frac{N}{p^2}<4$ for $p>5$. But then $4\nmid N$, a contradiction. So $p^2\nmid N$, i.e. the
$p$-Sylow subgroup of $\jac{C}(\FF_p)$ is of order at most $p$.
\end{proof}

Now consider the case $p\leq 5$. Assume at first that $D\equiv 2,3\mod{4}$. Since
$\omega_1\overline\omega_1=\omega_2\overline\omega_2=p$, we know that $|c_1\pm c_2\sqrt{D}|\leq\sqrt{p}$. Thus,
    \begin{align*}
    |c_2\sqrt{D}| &=    \frac{1}{2}\left|c_1+c_2\sqrt{D}-\left(c_1-c_2\sqrt{D}\right)\right| \\
           &\leq \frac{1}{2}\left(\left|c_1+c_2\sqrt{D}\right|+\left|c_1-c_2\sqrt{D}\right|\right) \\
       &\leq \sqrt{p}.
    \end{align*}
Similarly we see that $|c_1|\leq\sqrt{p}$. Assume that $D>5$. Then $|c_2|\leq\sqrt{\frac{p}{D}}<1$. So $c_2=0$, since
$c_2\in\ZZ$. This contradicts remark~\ref{rem:deg4}, i.e. $D\leq 5$. Now assume that $D=2$. Then
$c_2\leq\sqrt{\frac{p}{2}}\leq\sqrt{\frac{5}{2}}$, i.e. $c_2\in\{0,\pm 1\}$. Therefore it follows by calculating $N$
for each of the possible values of $c_1$ and $c_2$, that if $p^2\mid N$, then $c_2=0$. This is again a contradiction.
So if $D=2$, then $p^2\nmid N$. Similar it follows that if $D=3$, then $p^2\nmid N$.

Finally assume that $D\equiv 1\pmod{4}$. Then it follows from $\omega_1\overline\omega_1=\omega_2\overline\omega_2=p$
that $|c_1+c_2\frac{1\pm\sqrt{D}}{2}|\leq\sqrt{p}$. Thus, $|c_2\sqrt{D}|\leq 2\sqrt{p}$ and $|2c_1-c_2|\leq 2\sqrt{p}$.
Assume that $D>20$. Then $|c_2|<2\sqrt{\frac{5}{20}}=1$, i.e. $c_2=0$, a contradiction. So $D\leq 20$. By calculating
$N$ for each of the possible values of $p$, $D$, $c$ and $c_2$ it follows that $p^2\nmid N$ also in this case. Hence
the following lemma is established.

\begin{lem}\label{lem:2}
Let $C$ be a hyperelliptic curve of genus two defined over a prime field~$\FF_p$ of characteristic $p\leq 5$. Assume
that $\End(C)\simeq\heltal{K}$ and that the $p$-power Frobenius under this isomorphism is given by a number in
$\heltal{K_0}+\eta\heltal{K_0}$, where $\eta$ is given as in remark~\ref{rem:quarticCM}. Then the $p$-Sylow
subgroup of $\jac{C}(\FF_p)$ is either trivial or of order~$p$.
\end{lem}

Summing up, the following theorem holds.

\begin{teo}\label{teo}
Let $C$ be a hyperelliptic curve of genus two defined over a prime field~$\FF_p$. Assume that $\End(C)\simeq\heltal{K}$
and that the $p$-power Frobenius under this isomorphism is given by a number in $\heltal{K_0}+\eta\heltal{K_0}$,
where $\eta$ is given as in remark~\ref{rem:quarticCM}. Then the $p$-Sylow subgroup of $\jac{C}(\FF_p)$ is either
trivial or of order~$p$.
\end{teo}


\begin{thebibliography}{99}

    \bibitem[Atkin and Morain(1993)]{atkin-morain} \textsc{A.O.L. Atkin and F. Morain}. Elliptic curves and primality
    proving. \emph{Math. Comp.}, vol.~61, pp.~29--68, 1993.
    \bibitem[Cassels and Flynn(1996)]{cassels} \textsc{J.W.S. Cassels and E.V. Flynn}. \emph{Prolegomena to a Middlebrow Arithmetic of Curves of Genus $2$}. London Mathematical Society Lecture Note Series. Cambridge University Press, 1996.
    \bibitem[Gaudry \emph{et al}(2005)Pierrick Gaudry]{gaudry} \textsc{P. Gaudry, T. Houtmann, D. Kohel, C. Ritzenthaler and A. Weng}.
    The $p$-adic CM-Method for Genus $2$. 2005. \url{http://arxiv.org}.
    \bibitem[Koblitz(1989)]{koblitz89} \textsc{N. Koblitz}. Hyperelliptic cryptosystems. \emph{J. Cryptology}, vol.~1,
    pp.~139--150, 1989.
    \bibitem[Lang(1959) Serge Lang]{lang59} \textsc{S. Lang}. \emph{Abelian Varieties}. Interscience, 1959.
    \bibitem[Shimura(1998) Goro Shimura]{shi} \textsc{G. Shimura}. \emph{Abelian Varieties with Complex Multiplication and Modular Functions}.
    Princeton University Press, 1998.
    \bibitem[Spallek(1994)]{spallek} \textsc{A.-M. Spallek}. \emph{Kurven vom Geschlecht $2$ und ihre Anwendung in
    Public-Key-Kryptosystemen}. PhD thesis, Institut für Experimentelle Mathe\-matik, Universität GH Essen, 1994.
    \bibitem[Tate(1966)]{tate} \textsc{J. Tate}. Endomorphisms of abelian varieties over finite fields. \emph{Invent. Math.}, vol.~2, pp.~134--144, 1966.
    \bibitem[Weng(2003)]{weng03} \textsc{A. Weng}. Constructing hyperelliptic curves of genus~$2$ suitable for
    crypto\-graphy. \emph{Math. Comp.}, vol.~72, pp.~435--458, 2003.
\end{thebibliography}
\end{document}